\documentclass[runningheads]{lncse}

\usepackage{amsmath}
\usepackage{amsfonts}
\usepackage{epsfig}
\usepackage{graphics}
\usepackage{bm}
\usepackage{algorithm}
\usepackage{algorithmic}
\usepackage{url}
\usepackage{booktabs}
\usepackage{colortbl}

\begin{document}

\title{Block variants of the COCG and COCR methods for solving complex symmetric linear
systems with multiple right-hand sides}

\titlerunning{Block versions of COCG and COCR}

\author{Xian-Ming Gu\inst{1,2} \and Bruno Carpentieri\inst{3} \and Ting-Zhu Huang\inst{1}
\and Jing Meng\inst{4}}

\authorrunning{X.-M. Gu, B. Carpentieri, T.-Z. Huang, and J. Meng}

\institute{
School of Mathematical Sciences, University of Electronic Science and Technology of China, Chengdu
611731, P.R. China, {\tt guxianming@live.cn}, {\tt tingzhuhuang@126.com}
\and
Institute of Mathematics and Computing Science, University of Groningen,
Nijenborgh 9, P.O. Box 407, 9700 AK Groningen, The Netherlands
\and
School of Science and Technology, Nottingham Trent University, Clifton Campus, Nottingham, NG11 8NS,
United Kingdom, {\tt bcarpentieri@gmail.com}
\and
School of Mathematics and Statistics, Taishan University, Taian, 271021, P.R. China,
{\tt mengmeng-erni@163.com}
}

\maketitle

\begin{abstract}
In the present study, we establish two new block variants of the Conjugate Orthogonal 
Conjugate Gradient (COCG) and the Conjugate $A$-Orthogonal Conjugate Residual (COCR) 
Krylov subspace methods for solving complex symmetric linear systems with multiple 
right hand sides. The proposed Block iterative solvers can fully exploit the complex 
symmetry property of coefficient matrix of the linear system. We report on extensive 
numerical experiments to show the favourable convergence properties of our newly 
developed Block algorithms for solving realistic electromagnetic simulations.
\end{abstract}

\section{Introduction}
\label{author_mini4:sec:2}
In this paper we are interested in the efficient solution of linear systems
with multiple right-hand sides (RHSs) of the form
\begin{equation}
AX = B,\quad\ A\in \mathbb{C}^{n\times n},\ \ X,B\in \mathbb{C}^{n\times p},\ p\ll n,
\label{eq1.1}
\end{equation}
where $A$ is a non-Hermitian but symmetric matrix, i.e, $A \neq A^H$ and $A = A^T$. 
Linear systems of this form arise frequently in electromagnetic scattering applications, 
for example in monostatic radar cross-section calculation, where each right-hand side 
typically corresponds to an incident wave illuminating the target at a given angle 
of incidence~\cite{ISDLG,CDGS05}.

Roughly speaking, computational techniques for solving linear systems on modern 
computers can be divided into the class of direct and of iterative methods. Block 
iterative Krylov subspace methods are particularly designed for solving efficiently 
linear systems with multiple RHSs~(cf. \cite{JZJZA,MHGB}). Block algorithms require 
one or more matrix product operations of the form $AV$, with $V\in \mathbb{C}^{n\times 
p}$ an arbitrary rectangular matrix, per iteration step. Thus they can solve the 
typical memory bottlenecks of direct methods. However, most of them, such as the 
Block Bi-Conjugate Gradient (bl\_bicg)~\cite{DPOL}, Block Bi-Conjugate Residual 
(bl\_bicr)~\cite{JZJZA}, Block BiCGSTAB (bl\_bicgstab)~\cite{AGKJHS}, Block BiCRSTAB 
(bl\_bicrstab)~\cite{JZJZA}, Block QMR (bl\_qmr)~\cite{RWFMM}, Block IDR($s$) 
(bl\_idr($s$)) \cite{blidrs} and Block GMRES (bl\_gmres)~\cite{BVED} methods, do 
not naturally exploit any symmetry of $A$.

Methods that can exploit the symmetry of $A$ are typically of (quasi) minimal residual type 
(i.e. bl\_sqmr)~\cite{RWFMM}. Tadano and Sakurai recently proposed the Block COCG 
(bl\_cocg) \cite{HTTSA} method, which can be regarded as a natural extension of the 
COCG \cite{HAVVJBM} algorithm for solving linear systems (\ref{eq1.1}). Both these 
two methods need one operation $AV$ per iteration step. In this paper we revisit 
the Block COCG method, presenting a more systematic derivation than the one presented 
\cite{HTTSA}, and we introduce a new Block solver~(bl\_cocr) that can be seen as an 
extension of the COCR algorithm proposed in~\cite{TSSLZ}. The numerical stability 
of the bl\_cocg and the bl\_cocr methods are enhanced by the residual orthonormalization 
technique~\cite{AADR}.

The paper is organized as follows. In Section~\ref{author_mini4:sec:2x} we present the general framework for the development of the bl\_cocg and the bl\_cocr solvers. In Section 3 we study their numerical stability properties and then we show how  to improve their convergence by employing the residual orthonormalization technique. In Section 3, we report on extensive numerical experiments to illustrate the effectiveness of the two new iterative methods in computational electromagnetics. Finally, some conclusions arising from this work are presented in Section 4.

\section{The Block COCG and Block COCR methods}
\label{author_mini4:sec:2x}
Let $X^{m + 1}\in\mathbb{C}^{n\times p}$ be the $(m + 1)$th approximate solution of linear systems (\ref{eq1.1}) satisfying the following condition
\begin{equation}
X_{m + 1} = X_0 + Z_{m + 1},\quad Z_{m +1}\in \mathcal{K}^{\diamond}_{m+1}(A; R_0),
\label{eq1.2}
\end{equation}
where $R_0 = B - AX_0$ is an initial residual and $\mathcal{K}^{\diamond}_{m+1}(A; R_0)$ is the block
Krylov subspace~\cite{MHGB} defined as
\begin{equation}
\mathcal{K}^{\diamond}_{m+1}(A; R_0) = \Big\{\sum^{m}_{j = 0}A^j R_0\gamma_j\mid \gamma_j \in
\mathbb{C}^{p\times p}~(j = 0,1,\ldots,m)\Big\}.
\label{eq1.3}
\end{equation}
Compared with conventional Krylov subspace methods, where ${\bm x}^{(j)}_{m + 1} - {\bm x}^{(j)}_0
\in\mathcal{K}_{m + 1}(A, {\bm r}^{(j)}_0)$, note that block Krylov methods can
search the approximate solutions into larger spaces, and thus they may require less
iterations to converge to a given accuracy. In the next section we introduce the framework for the development of the Block COCG and the Block COCR methods.

\subsection{Derivation of the Block COCG and Block COCR methods}
\label{author_mini4:subsec:2}
According to Eqs. (\ref{eq1.2})--(\ref{eq1.3}), the $(m + 1)$th residual $R_{m+1} =
B - AX_{m +1}$ of the Block COCG method \cite{HTTSA} and the Block COCR method is
computed by the following recurrence relations,
\begin{eqnarray}
R_0 = P_0 = B - AX_0 \in \mathcal{K}^{\diamond}_1(A; R_0), \nonumber\\
R_{m+1} = R_m - AP_m\alpha_m\in \mathcal{K}^{\diamond}_{m + 2}(A; R_0),\nonumber\\
P_{m + 1} = R_{m + 1} + P_m \beta_m \in \mathcal{K}^{\diamond}_{m + 2}(A;
R_0).
\label{author_mini4:eq:01x}
\end{eqnarray}
Here, $P_{m+1} \in \mathbb{C}^{n\times p}, \alpha_m, \beta_m\in \mathbb{C}^{p\times
p}$. The $(m + 1)$th approximate solution $X_{m+1}$ is updated through the recurrence
relation
\begin{equation}
X_{m+1} = X_m + P_m \alpha_m.
\label{author_mini4:eq:01xx}
\end{equation}
Similarly to the framework introduced in \cite{XMGMC}, different formulae for the
$p\times p$ matrices $\alpha_m,\beta_m~(m = 0,1,\ldots)$ in the recurrences
(\ref{author_mini4:eq:01x})--(\ref{author_mini4:eq:01xx}) lead to different iterative
algorithms. Denoting by $\mathcal{L}$ the \textit{block constraints subspace}, these matrices
$\alpha_m,\beta_m$ are determined by imposing the orthogonality conditions
\begin{equation}
R_m \perp \mathcal{L}\quad\ \mathrm{and}\quad\ AP_m \perp \mathcal{L}.
\label{author_mini4:eq:01y}
\end{equation}

The Block COCG and the Block COCR methods correspond to the choices $\mathcal{L}=\mathcal{K}^{\diamond}_m(\bar{A}; \bar{R}_0)$
and $\mathcal{L}=\bar{A}\mathcal{K}^{\diamond}_m(\bar{A};\bar{R}_0)$, respectively. In Table~\ref{author_mini4:tab:1}, the
conjugate orthogonality conditions imposed to determine $\alpha_m$ and $\beta_m$ are summarized for the sake of clarity.
\begin{table}
\centering
\caption{Orthogonality conditions imposed to determine $p\times p$ matrices $\alpha_m,\beta_m$}
\label{author_mini4:tab:1}
\begin{tabular}{p{2cm}p{3cm}p{3cm}}
\hline\noalign{\smallskip}
Matrix             & Block COCG                             & Blcok COCR   \\
$\alpha_m,\beta_m$ & $R_m \perp \mathcal{K}^{\diamond}_m(\bar{A};\bar{R}_0)$  & $R_m
\perp \bar{A}\mathcal{K}^{\diamond}_m(\bar{A};\bar{R}_0)$ \\
                   & $AP_m \perp \mathcal{K}^{\diamond}_m(\bar{A};\bar{R}_0)$ & $AP_m
                   \perp \bar{A}\mathcal{K}^{\diamond}_m(\bar{A};\bar{R}_0)$ \\
\noalign{\smallskip}\hline\noalign{\smallskip}
\end{tabular}
\end{table}

We show the complete Block COCR algorithm in Algorithm \ref{alg1x}. We use
the notation $\|\cdot\|_F$ for the Frobenius norm of a matrix, and $\epsilon$ is a sufficiently
small user-defined value. We see that the Block COCR method requires two matrix products
$AP_{m + 1}$, $AR_{m + 1}$ at each iteration step. While the product $AR_{m + 1}$ is computed by explicit matrix multiplication, the product $AP_{m + 1}$ is computed by the recurrence relation at line 9, to reduce the computational complexity.
Note that the Block COCG and the Block COCR methods can be derived from the Block BiCG and the Block BiCR methods, respectively, by choosing the initial auxiliary residual $\hat{R}_0 = \bar{R}_0$ and removing some redundant computations; we refer to the recent work~\cite{XMGMC} for similar discussions about the derivation of conventional non-block Krylov subspace methods for complex symmetric linear systems with single RHS.

\begin{algorithm}
\caption{The Block COCR method}
\begin{algorithmic}[1]
  \STATE $X_0\in \mathbb{C}^{n\times p}$ is an initial guess, $R_0 = B - AX_0$,
  \STATE Set $P_0 = R_0$, $U_0 = V_0 = AR_0$,
  \FOR{$m = 0,1,\ldots$, until $\|R_m\|_F/\|R_0\|_F \leq \epsilon$}
  \STATE Solve $(U^{T}_mU_m)\alpha_m = R^{T}_m V_m$ for $\alpha_m$,
  \STATE $X_{m + 1} = X_m + P_m \alpha_m$,
  \STATE $R_{m + 1} = R_m - U_m \alpha_m$ and $V_{m+1} = AR_{m+1}$,
  \STATE Solve $(R^{T}_mV_m)\beta_m = R^{T}_{m + 1} V_{m + 1}$ for $\beta_m$,
  \STATE $P_{m + 1} = R_{m + 1} + P_m \beta_m$,
  \STATE $U_{m + 1} = V_{m + 1} + U_m \beta_m$,
  \ENDFOR
\end{algorithmic}
\label{alg1x}
\end{algorithm}

\subsection{Improving the numerical stability of the Block COCG and Block COCR methods by residual orthonormalization}
One known problem with Block Krylov subspace methods is that the residual norms  may not converge when the number $p$ of right-hand sides
is large, mainly due to numerical instabilities, see e.g.~\cite{AADR}. These instabilities often arise because of the loss of linear independence among the column vectors of the $n\times p$ matrices that appear in the methods, such as $R_{m}$ and $P_m$. Motivated by this concern, in this section we propose to use the residual orthonormalization technique to enhance the numerical stability
of the Block COCG and Block COCR algorithms. This efficient technique was introduced in~\cite{AADR} in the context of the Block CG method \cite{DPOL}.

Let the Block residual $R_m$ be factored as $R_m = Q_m\xi_m$ by conventional QR
factorization\footnote{For our practical implementation, we use MATLAB qr-function
``\texttt{qr}($W$,0)" for a given matrix $W\in\mathbb{C}^{n\times p}$.}, with $Q^{H}_m Q_m = I_p$. Here $I_p$ denotes the identity
matrix of order $p$ and $\xi_m\in \mathbb{C}^{p\times p}$. From (\ref{author_mini4:eq:01x}), the following equation can be obtained
\begin{equation}
Q_{m +1}\tau_{m + 1} = Q_m - AS_m \alpha'_k.
\end{equation}
Here, $\tau_{m +1} \equiv \xi_{m + 1}\xi_{m - 1}$, $\alpha'_k \equiv \xi_m \alpha_m
\xi_{m - 1}$, and $S_m = P_m\xi_{m - 1}$. In the new Algorithms~\ref{alg2x}-\ref{alg3x}, 
the matrix $\beta'_m$ is defined as $\alpha'_m\equiv \xi_m \beta_m \xi^{-1}_{m + 1}$. 
The residual norm is monitored by $\|\xi_m\|_F$ instead of $\|R_m\|_F$, since the 
Frobenius norm of $R_m$ satisfies $\|R_m\|_F = \|\xi_m\|_F$. Note that the QR decomposition 
is calculated at each iteration. However, the numerical results shown in the next 
section indicate that the extra cost is amortized by the improved robustness of the 
two Block solvers.

\begin{algorithm}[t]
\caption{Algorithm of the Block COCG method with residual orthonormalization (bl\_cocg\_rq)}
\begin{algorithmic}[1]
  \STATE $X_0\in \mathbb{C}^{n\times p}$ is an initial guess, and compute $Q_0\xi_0 = B - AX_0$,
  \STATE Set $S_0 = Q_0$,
  \FOR{$m = 0,1,\ldots$, until $\|\xi_m\|_F/\|B\|_F \leq \epsilon$}
  \STATE Solve $(S^{T}_m A S_m)\alpha'_m = Q^{T}_m Q_m$ for $\alpha'_m$,
  \STATE $X_{m + 1} = X_m + S_m \alpha'_m\xi_m$,
  \STATE $Q_{m + 1}\tau_{m + 1} = Q_m - AS_m \alpha'_m$ and $\xi_{m+1} = \tau_{m+1}\xi_m$,
  \STATE Solve $(Q^{T}_mQ_m)\beta'_m = \tau^{T}_{m + 1} Q^{T}_{m + 1}Q_{m + 1}$ for $\beta'_m$,
  \STATE $S_{m + 1} = Q_{m + 1} + S_m \beta'_m$,
  \ENDFOR
\end{algorithmic}
\label{alg2x}
\end{algorithm}

\begin{algorithm}
\caption{Algorithm of the Block COCR method with residual orthonormalization (bl\_cocr\_rq)}
\begin{algorithmic}[1]
  \STATE $X_0\in \mathbb{C}^{n\times p}$ is an initial guess, and compute $Q_0\xi_0 = B - AX_0$,
  \STATE Set $S_0 = Q_0$ and $U_0 = V_0 =  AQ_0$,
  \FOR{$m = 0,1,\ldots$, until $\|\xi_m\|_F/\|B\|_F \leq \epsilon$}
  \STATE Solve $(U^{T}_mU_m)\alpha'_m = Q^{T}_m U_m$ for $\alpha'_m$,
  \STATE $X_{m + 1} = X_m + P_m \alpha'_m$
  \STATE $Q_{m + 1}\tau_{m + 1} = Q_m - U_m \alpha'_m$ and $\xi_{m+1} = \tau_{m+1}\xi_m$,
  \STATE Compute $V_{m + 1} = AQ_{m + 1}$,
  \STATE Solve $(Q^{T}_mV_m)\beta_m = \tau^{T}_{m + 1} Q^{T}_{m + 1}V_{m + 1}$ for $\beta'_m$,
  \STATE $S_{m + 1} = Q_{m + 1} + S_m \beta'_m$,
  \STATE $U_{m + 1} = V_{m + 1} + U_m \beta'_m$,
  \ENDFOR
\end{algorithmic}
\label{alg3x}
\end{algorithm}

\section{Numerical experiments}
\label{author_mini4:subsec:3}
In this section, we carry out some numerical experiments to show the potential
effectiveness of the proposed iterative solution strategies in computational electromagnetics. We compare the bl\_cocg,
bl\_cocg\_rq, bl\_cocr, bl\_cocr\_rq methods against other popular block Krylov subspace
methods such as bl\_qmr, bl\_bicgstab, bl\_bicrstab, bl\_idr($s$) (selecting matrix $P =
rand(n,sp)$, see~\cite{blidrs}) and restarted bl\_gmres(m).
We use the value $m = 80$ for the restart in bl\_gmres(m). The experiments have
been carried out in double precision floating point arithmetic with MATLAB 2014a
(64 bit) on PC-Intel(R) Core(TM) i5-3470 CPU 3.20 GHz, 8 GB of RAM.

The different Block algorithms are compared in terms of number of iterations, denoted as
\textit{Iters} in the tables, and $\log_{10}$ of the final true relative residual
norm defined as $\log_{10}(\|B - AX_{\mathrm{final}}\|_F/\|B\|_F)$, denoted as
\textit{TRR}. The iterative solution is started choosing $X_0 = O\in\mathbb{C}^{n\times p}$ 
as initial guess. The stopping criterion in our runs is the reduction of the norm 
of the initial Block residual by eight orders of magnitude, i.e., $\|R_m\|_F/\|B\|_F 
\leq Tol = 10^{-10}$. The right-hand side $B$ is computed by the MATLAB function 
{\texttt{rand}}. In the tables, the symbol ``$\dag$" indicates no convergence 
within $n$ iterations, or $n/m$ cycles for the bl\_gmres($m$) method.

The first test problems are three matrices extracted from the Matrix Market
collection\footnote{\url{http://math.nist.gov/MatrixMarket/matrices.html}},
arising from modeling acoustic scattering problems. They are denoted as~\texttt{young1c, young2c}, 
and {\tt young3c}. The results of our experiments are presented in Table~\ref{tab2}. 
The symbol ${}^{*}$ used for the bl\_bicgstab, bl\_idr(4), and bl\_bicrstav methods 
indicate that these three methods require no less than two matrix products $AV$ per 
iteration step. The symbol ${}^{**}$ refers to the number of outer iterations in the 
Block GMRES($m$) method, when it can achieve convergence; refer to \cite{HXZGW} for details. This notation is used throughout this section.

\begin{table}[t]\footnotesize\tabcolsep=4.3pt
\begin{center}
\caption{The numerical results of different iterative solvers for the first example.}
\begin{tabular}{cccccccccc}
\hline Method &\multicolumn{3}{c}{\texttt{young2c} ($p = 10$)}&\multicolumn{3}{c}{
\texttt{young3c} ($p = 8$)}&\multicolumn{3}{c}{\texttt{young1c} ($p = 8$)}\\
[-2pt]\cmidrule(r{0.5em}){2-4} \cmidrule(l{0.5em}r{0.5em}){5-7}\cmidrule(l{0.5em}){8-10}
\\[-11pt]
            &$Iters$     &{\it TRR} &CPU   &$Iters$ &{\it TRR} &CPU    &$Iters$   &{\it TRR} &CPU    \\
\hline
bl\_cocg      &238       &-10.03    &0.17  &$\dag$  &$\dag$    &$\dag$ &329       &-10.16    &0.16 \\
bl\_cocg\_rq  &142       &-10.14    &0.13  &151     &-10.00    &0.09   &177       &-10.29    &0.12 \\
bl\_cocr      &201       &-10.07    &0.15  &145     &-9.95     &0.04   &221       &-10.07    &0.12 \\
bl\_cocr\_rq  &138       &-10.18    &0.13  &146     &-10.03    &0.05   &180       &-10.18    &0.13 \\
bl\_sqmr      &154       &-9.87     &0.29  &131     &-10.39    &0.09   &188       &-9.88     &0.25 \\
bl\_bicgstab  &395$^{*}$ &-10.09    &0.41  &$\dag$  &$\dag$    &$\dag$ &433$^{*}$ &-10.04    &0.35 \\
bl\_bicrstab  &356$^{*}$ &-9.96     &0.46  &$\dag$  &$\dag$    &$\dag$ &417$^{*}$ &-9.71     &0.44 \\
bl\_idr(4)    &269$^{*}$ &-8.57     &0.28  &$\dag$  &$\dag$    &$\dag$ &334$^{*}$ &-10.10    &0.27 \\
bl\_gmres(m)  &3$^{**}$  &-10.08    &24.5  &$\dag$  &$\dag$    &$\dag$ &$\dag$    &$\dag$    &$\dag$ \\
\hline
\end{tabular}
\label{tab2}
\end{center}
\end{table}

Table \ref{tab2} shows the results with nine different Block Krylov solvers. Although the bl\_cocg and bl\_cocr methods required more $Iters$, they are
more competitive than the bl\_sqmr method in terms of CPU time and \textit{TRR} (except the case of \texttt{young3c}). Bl\_cocr method is more robust than  bl\_cocg in terms of $Iters$, CPU time and \textit{TRR}. The bl\_cocg\_rq and bl\_cocr\_rq variant are very efficient in terms of \textit{TRR} and CPU time. The bl\_bicgstab, bl\_bicrstab, bl\_idr(4), and bl\_gmres($m$) methods cannot solve the test problem (\texttt{young3c}), while bl\_cocg and bl\_cocr converge rapidly. Due to the long iterative recurrence, the bl\_gmres($m$) method is typically expensive.

In the second experiment we consider three dense matrices arising from monostatic radar cross-section calculation; they are denoted as \texttt{sphere2430}, \texttt{parallelepipede}, {\texttt{cube1800}}. These problems are available from our GitHub repository\footnote{\url{https://github.com/Hsien-Ming-Ku/Test_matrices/tree/master/Example2}}, and we choose $p = 8$. Although rather small, the selected dense problems are representative of realistic radar-cross-section calculation~\cite{CDGS05}. Larger problems would require a Fortran or C implementation of the solvers and will be considered in a separate study.
Numerical results for each test problem are summarized in Table~\ref{tab3}.

\begin{table}[t]\footnotesize\tabcolsep=4.3pt
\begin{center}
\caption{The numerical results of different iterative solvers for Example 1.}
\begin{tabular}{cccccccccc}
\hline Method &\multicolumn{3}{c}{\texttt{sphere2430}}&\multicolumn{3}{c}{\texttt{parallelepipede}}
&\multicolumn{3}{c}{\texttt{cube1800}}\\
[-2pt]\cmidrule(r{0.5em}){2-4} \cmidrule(l{0.5em}r{0.5em}){5-7}\cmidrule(l{0.5em}){8-10}
\\[-11pt]
            &$Iters$     &{\it TRR} &CPU   &$Iters$   &{\it TRR} &CPU    &$Iters$   &{\it TRR} &CPU    \\
\hline
bl\_cocg      &189       &-10.07    &4.16  &176       &-10.02    &2.40   &174       &-10.21    &1.94 \\
bl\_cocg\_rq  &169       &-10.00    &3.77  &156       &-10.13    &2.13   &156       &-10.08    &1.74 \\
bl\_cocr      &186       &-10.03    &4.12  &174       &-10.02    &2.35   &169       &-10.00    &1.84 \\
bl\_cocr\_rq  &166       &-10.05    &3.77  &152       &-10.15    &2.11   &151       &-10.09    &1.73 \\
bl\_sqmr      &172       &-9.84     &4.15  &161       &-9.91     &2.42   &159       &-9.97     &2.11 \\
bl\_bicgstab  &379$^{*}$ &-10.04    &16.5  &370$^{*}$ &-10.04    &9.94   &396$^{*}$ &-10.29    &8.42 \\
bl\_bicrstab  &392$^{*}$ &-9.57     &17.3  &355$^{*}$ &-9.85     &9.98   &303$^{*}$ &-8.38     &6.70 \\
bl\_idr(4)    &409$^{*}$ &-9.64     &22.1  &474$^{*}$ &-10.11    &16.5   &334$^{*}$ &-9.43     &10.2 \\
bl\_gmres(m)  &2$^{**}$  &-10.07    &38.2  &2$^{**}$  &-10.04    &33.3   &2$^{**}$  &-10.09    &22.1  \\
\hline
\end{tabular}
\label{tab3}
\end{center}
\end{table}

Table \ref{tab3} displays the results with again nine different Block Krylov solvers. We can see that the bl\_sqmr method requires less \textit{Iters} to converge compared to the bl\_cocg and bl\_cocr methods. However, it is more expensive in terms of CPU time except on the \texttt{sphere2430} problem. Besides, the true residual norms produced by the bl\_sqmr method are larger than those of both bl\_cocg and bl\_cocr. Furthermore, bl\_cocg\_rq and bl\_cocr\_rq are the most effective and promising solvers in terms of \textit{Iters} and CPU time. Specifically, the bl\_cocr\_rq method is slightly more efficient than the bl\_cocg\_rq method in terms of \textit{TRR}.

\section{Conclusions}
\label{author_mini4:subsec:4}
In this paper, a framework for constructing new Block iterative Krylov subspace 
methods is presented. Two new matrix solvers that can exploit the symmetry of $A$ 
for solving complex symmetric non-Hermitian linear systems~(\ref{eq1.1}) are introduced. 
Stabilization techniques based on residual orthonormalization strategy are discussed 
for both methods. The numerical experiments show that the solvers can be viable 
alternative to standard Krylov subspace methods for solving  complex symmetric linear systems 
with multiple RHSs efficiently. Obviously, for solving realistic electromagnetic 
problems they both need to be combinated with suitable preconditioners that reflect 
the symmetry of $A$; we refer the reader to, e.g.,~\cite{CABO12,PLRRSC,CDGM01B} 
for some related studies.

\ifx\undefined\bysame
\newcommand{\bysame}{\leavevmode\hbox to3em{\hrulefill}\,}
\fi

\end{document}